\documentclass[journal]{IEEEtran}

\ifCLASSINFOpdf
\else
\fi
\usepackage{lipsum}
\usepackage{amsfonts}
\usepackage{amsmath}
\usepackage{amsthm}
\usepackage{amssymb}
\usepackage{graphicx}
\usepackage{epstopdf}
\usepackage{comment}
\usepackage{algorithmic}
\usepackage{color}
\usepackage[]{algorithm2e,setspace}
\ifpdf
  \DeclareGraphicsExtensions{.eps,.pdf,.png,.jpg}
\else
  \DeclareGraphicsExtensions{.eps}
\fi

\newcommand{\TheTitle}{A multilevel based reweighting algorithm with joint regularizers for sparse recovery}

\title{{\TheTitle}}

\author{Jackie Ma \thanks{J. Ma is with the   Image   and   Video   Coding
Group,   Fraunhofer   Institute   for   Telecommunications--Heinrich   Hertz   Institute,  Berlin  10587,  Germany  (e-mail:  jackie.ma@hhi.fraunhofer.de).}
\quad 
Maximilian M\"{a}rz \thanks{M. M\"{a}rz is with the Department
of Mathematics,   Technische Universit\"{a}t Berlin, Berlin, 10623 Berlin, Germany 
(e-mail: maerz@math.tu-berlin.de).

This work has been submitted to the IEEE for possible publication.
Copyright may be transferred without notice, after which this version may
no longer be accessible}}

\DeclareMathOperator{\diag}{diag}

\usepackage{float}
\usepackage{graphicx}
\usepackage{tabularx}

\newtheoremstyle{j}%
{3pt}%
{3pt}%
{}%
{\parindent}%
{\bfseries}%
{.}%
{.5em}%
{}%

\theoremstyle{plain}

\newtheorem*{rem*}{Remark}
\newtheorem*{concl*}{Conclusion}
\newtheorem*{theorem*}{Theorem}
\newtheorem*{cor*}{Corollary}
\newtheorem*{algo*}{Algorithm}

\newtheorem{theorem}{Theorem}[section]

\newtheorem{deff}[theorem]{Definition}

\newcommand{\Z}{\mathbb{Z}}
\newcommand{\N}{\mathbb{N}}

\newcommand{\R}{\mathbb{R}}

\newcommand{\C}{\mathbb{C}}

\newcommand{\psitilde}{\widetilde{\psi}}

\newcommand{\eps}{\varepsilon}

\newcommand{\bal}{\begin{align}}
\newcommand{\eal}{\end{align}}

\newcommand{\bM}{\begin{pmatrix}}
\newcommand{\eM}{\end{pmatrix}}

\newcommand{\norm}[1]{\left\lVert#1\right\rVert}

\DeclareMathOperator*{\argmin}{arg\min}

\DeclareMathOperator{\TGV}{TGV}
\DeclareMathOperator{\divergenz}{div}
\DeclareMathOperator{\cur}{cur}
\DeclareMathOperator{\maxIter}{maxIter}

\begin{document}

\maketitle

\begin{abstract}

We propose an algorithmic framework based on ADMM/split Bregman that combines a \emph{multilevel adapted, iterative reweighting strategy} and a second \emph{total generalized variation} regularizer. The level adapted reweighting strategy is a combination of reweighted $\ell^1$-minimization and additional compensation factors for a uniform treatment of the sparsity structure across all levels. Classical multilscale transforms that are very well suited for this algorithm are, for instance, the \emph{wavelet transform} and the \emph{shearlet transform}. The proposed algorithm is tested for the reconstruction of images from their Fourier measurements and Radon measurements, respectively. The numerical experiments show a highly improved performance at relatively low additional computational costs compared to many other well established methods. 
\end{abstract}

\section{Introduction}
The field of sparse recovery has a wide range of applications in many different areas such as medical imaging \cite{LusDonPau,CTRef}, astronomy \cite{Astronomy1,Astronomy2}, electron microscopy \cite{EM1} etc. One of the great successes in this area are new developments of multiscale sparsifying transforms since the invention of wavelets and the wavelet transform. Indeed wavelets are known to compress natural images very effectively since most natural images are sparse in a wavelet domain. However, they lack in directional sensitivity and are therefore not optimal for images that are governed by curvilinear structures. For precisely that reason, almost a decade ago directional systems such as curvelets \cite{CanDon} and shearlets \cite{KitKutLim, GuoKutLab2006} have been created to overcome this deficit. We will briefly recall the concepts of wavelets and shearlets in Section \ref{sec:Section2} as both systems will play a significant role in the upcoming content of this paper.

The field of sparse recovery is predominantly influenced by the development of  \emph{compressed sensing} \cite{CanRomTao, Don}, a theory that guarantees the recovery of \emph{sparse signals} from incomplete measurements under the assumptions of \emph{sparsity} and \emph{incoherence}. These sparse signals are typically obtained by solving a convex optimization problem of the form
\begin{align}
    \min_u \| \Psi  u \|_1 \quad \text{ subject to } \quad Au = y, \label{eq:AnaMin}
\end{align}
where $u$ is the object of interest, $\Psi$ is a sparsifying transform, $A$ is a matrix representing the measurement process, and $y$ are the resulting measurements. If $u$ is already \emph{sparse} it suffices to let $\Psi$ to be the identity. Otherwise \eqref{eq:AnaMin} is called \emph{Basis Pursuit in the analysis formulation}.
 Very recently, Ahmad and Schniter considered a generalized variant of \eqref{eq:AnaMin} in \cite{AhmSch}, where they have used not only one sparsifying transform, but a composition of several sparsifying transforms. Furthermore, in order to improve the reconstruction quality the authors have combined their ideas with reweighted $\ell^1$-minimization into their framework. Indeed, the concept of \emph{reweighted $\ell^1$-minimization} introduced by Cand\`{e}s et al. in \cite{CanWakBoy} further promotes the sparsity of the recovered signal by iteratively updating a weighting matrix in the minimization problem. More precisely, one solves 
 \begin{align}
\min_u \| W_k\Psi u \|_1 \quad \text{ subject to } \quad Au = y, \label{eq:AnaMinRe}
\end{align}
iteratively for $k=1,2,\dots$ where after each iteration $k$ the diagonal weighting matrix $W_k$ is updated according to the sparsity structure of the current solution $u_{k+1}$ of \eqref{eq:AnaMinRe}. The role of $W_k$ is to mimic the actual sparsity structure of the true signal that one wishes to recover,. We will recap the ideas of reweighted $\ell^1$ later in Section \ref{sec:Section2} in more detail as this will also be one of the main ingredients of the algorithm that we propose in this work.

In order to find solutions of \eqref{eq:AnaMin} or approximations of such solutions, 
many different possibilities available in the literature. In this work, we focus on the \emph{split Bregman} algorithm \cite{SB,SB2}, which besides minor differences, is a reinvention of the \emph{alternating direction method of multipliers} (ADMM) \cite{GabMer, EckBer, BoyADMM}. It can also be seen as \emph{Douglas Rachford splitting} of the dual problem. For details on the relationship of these algorithms we refer the interested reader to \cite{Setzer,Esser}.  The split Bregman algorithm transforms the constrained problem \eqref{eq:AnaMin} into an unconstrained formulation and by introducing splitting variables they break the original problem down into, hopefully, easier ones. A key ingredient is then the so-called \emph{soft-thresholding}, \emph{shrinking} or \emph{proximal mapping} which gives a closed-form solution to some of the subproblems. One of the great advantages of ADMM and in particular split Bregman is that they are easily derived and very flexible in terms of  multiple regularizers. Furthermore, as we shall explain in this paper, it can be greatly combined with the idea of reweighted $\ell^1$. 
 
 It is often beneficial to consider a second sparsity promoting regularizer such as total variation (TV) \cite{SB,RecPF} to reduce artifacts coming from the sparsifying transform $\Psi$ in the 
solutions of sparse imaging problems  that are obtained via  \eqref{eq:AnaMin}. TV was initially proposed for denoising problems \cite{RudOshFat} and is very well established in image processing by now. However, it has been noticed that severe staircasing artifacts may  appear in images recovered by TV regularized reconstructions for large noise levels. In order to overcome this issue the authors of \cite{BreKunPoc} have introduced \emph{total generalized variation (TGV)} which is a generalization of TV to higher order derivatives. Since then TGV has been used in many applications \cite{BloTobUecFra,KnoBrePocRud} as a regularizer for inverse problems. 

 The consideration of a joint regularization scheme using shearlets and TGV has already been done by Guo et al. in \cite{GuoQinYin} in order to solve problems of the form \eqref{eq:AnaMin}. Our main contribution is to further exploit the general structure of the (multilevel) sparsity by combining iterative reweighting and adaptive multilevel weights associated to any multiscale transform, not only shearlets. This is also related to the approach in \cite{AhmSch}, but the algorithm and the conclusions derived in this work are different. In order to solve the resulting multilevel reweighted $\ell^1$ problem we make use of the flexibility of split Bregman by directly incorporating adaptive multilevel thresholds into the subproblem that is solved by a simple thresholding step. This results in an algorithm that comes with very little additional computational cost and almost automatically chosen regularization parameters. In particular, as we will show in this paper, the results are greatly improved compared to the non-reweighted analogue. We wish to mention, that the reweighting approach that we are considering in this paper is not to be confused with weighted $\ell^1$-minimization. Weighted $\ell^1$-minimization requires addition a priori knowledge to carefully design effetive weights. whereas in reweighting the weights are adaptive and do not require a priori knowledge. Furthermore, reweighted $\ell^1$-minimization is an iterative scheme and the minimization problem has to be solved several times (with updated weights) which is not the case for weighted $\ell^1$.

\begin{figure*}
 \includegraphics[width=.24\textwidth]{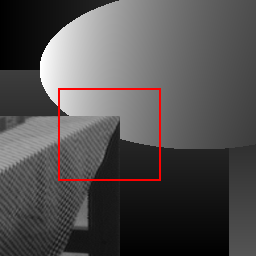}
 \includegraphics[width=.24\textwidth]{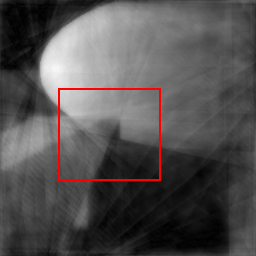}
 \includegraphics[width=.24\textwidth]{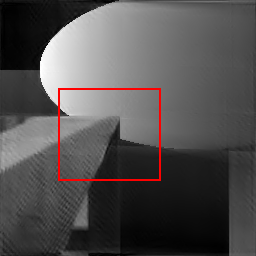}
 \includegraphics[width=.24\textwidth]{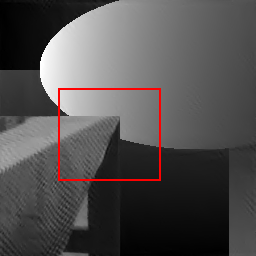}
\\[.5ex] 
 \includegraphics[width=.24\textwidth]{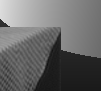}
 \includegraphics[width=.24\textwidth]{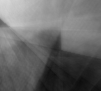}
 \includegraphics[width=.24\textwidth]{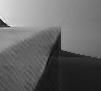}
 \includegraphics[width=.24\textwidth]{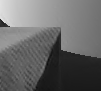}
 \caption{Reconstructions from 25 radial lines ($10.28\%$) through the k-space origin. \textbf{First column}: Original image with zoom. \textbf{Second column}: Reconstruction with inverse Fourier transform. Relative error: 0.146. Structured similarity index: 0.735. \textbf{Third column}: Reconstruction with redundant Daubechies 4 wavelets and total generalized variation regularizer without reweighting. Relative error: 0.060. Structured similarity index: 0.896. \textbf{Fourth column}: Reconstruction with the proposed method. Relative error: 0.031. Structured similarity index: 0.951.}\label{fig:abstract}
 \end{figure*}

 In Figure \ref{fig:abstract} we give a motivation and a first glance for the possible benefit of combining reweighting methods with, in this particular case, the wavelet transform. The reconstructions shown in Figure \ref{fig:abstract}are obtained from partial Fourier measurements of an synthetic test image which has parts that are certainly sparse in a wavelet dictionary and parts that are very well suited for TGV. In Section \ref{sec:Section4} we will present all details of the numerical implementation. Note that in Figure \ref{fig:abstract} the proposed method reduces the artifacts while still being able to reconstruct fine details using only $10.28\%$ of Fourier measurements.
 
\subsection{Outline}

In Section \ref{sec:Section2} we will give a compact overview of all methods that are needed in order to derive our proposed algorithm. In Section \ref{sec:Section3} we then present our ideas and the final method. The last section, Section \ref{sec:Section4} contains carefully conducted experiments with method, including important comparisons to other classical and novel algorithms that are known to work well for the recovery from incomplete Fourier measurements. Furthermore we will demonstrate the successful performance of our algorithm in the context of computed tomography.

\section{Sparse recovery, convex optimization, and sparsifying transforms}\label{sec:Section2}

In this section we present a short overview of current concepts and methods that are standard in the area of sparse recovery and are necessary to follow the rest of this work. For more details we refer the interested reader to the indicated literature.

\subsection{Compressed sensing and reweighted $\ell^1$-minimization}

The problem considered in compressed sensing can be explained by solving a system of underdetermined linear equations using prior information. Indeed, one is interested in solutions $u$ of the equation
\begin{align}
 Au = y, \label{eq:Axy}
\end{align}
where $y$ is a vector representing the acquired data, $A$ is a sensing matrix, and $u$ is the object of interest. The assumption that makes this problem in particular interesting is that $y$ should be of very 
small dimension compared to $u$ that lives in a much higher dimensional space, i.e., $A \in \C^{m \times n}$ with $m \ll n$. Furthermore, the aforementioned prior information that enables us to solve such an underdetermined 
system is \emph{sparsity}, i.e., although $u$ might be drawn from a much higher dimensional space only very few of its entries are nonzero. More precisely, a signal $u \in \C^n$ is called \emph{$s$-sparse} if
\begin{align*}
\| u \|_0 := \#\{ i \in \N \, : \, u(i) \neq 0, 1\leq i \leq n\}\leq s.
\end{align*}

A common approach to obtain sparse solutions of \eqref{eq:Axy} is to solve the following constrained convex optimization problem
\begin{align}
\label{min:standardl1}
  \min_u \| u\|_1 \quad \text{ subject to } \quad Au = y, 
\end{align}
see for example \cite{CanRomTao,FouRau}.

One of the possibilities to improve the recovery model \label{min:l1minimization} is to strengthen the effect of sparsity in the minimization problem. This can be done, for instance, by using the idea of reweighted $\ell^1$ introduced by Cand\`{e}s et al. in \cite{CanWakBoy} which can be described as follows. 

Suppose we are given measurements $y = Au \in \C^m$ of an $s$-sparse signal $$u = (u(1),\dots,u(n))^T \in \C^n$$ for a measurement matrix $A \in \C^{m \times n}$ for $m \ll n$.
When solving the minimization problem \eqref{min:standardl1} iteratively, one would ask for the following effect: large coefficients should be quickly identified and hence become ``cheaper'' in the minimization of the  objective function in \eqref{min:standardl1}, whereas very small coefficients 
should be neglected in the minimization since they are most likely going to be zero in the true signal. More precisely, let $u_0$ be the true signal and define a  diagonal weighting matrix  $W$ by
\begin{align}
 W(i,i) = \begin{cases}
                                    \frac{1}{|u_{0}(i)|},&  u_{0}(i) \neq 0 \\
                                    \infty, & u_{0}(i) = 0.
                                   \end{cases} \label{weights}
\end{align}
Now, if the signal $u_0$ was $s$-sparse, then under some assumptions \cite{CanWakBoy} the \emph{weighted} $\ell^1$-minimization problem
\begin{align*}
 \min_u \| W u\|_1 \quad\text{ subject to } \quad Au = y,
\end{align*}
will find the exact solution. However, since $u_0$ is usually unknown such weights are practically infeasible. Therefore Cand\`{e}s et al. proposed adaptive weights in \cite{CanWakBoy} that change at each iteration depending on the previously computed solution $u^k$ which is an approximation to $u_0$. More precisely, the following sequence of minimization problems are then considered
\begin{align*}
u_{k+1} =  \argmin_u \| W_k u\|_1 \quad \text{ subject to } \quad Au = y, 
\end{align*}
with a weighting matrix
\begin{align*}
 W_{k}(i,i) =   \frac{1}{|u_{k}(i)|+\varepsilon},
\end{align*}
where $\varepsilon > 0$ is a stability parameter and the initial weighting matrix $W^0$ is set to be the identity. In a series of numerical experiments it was shown in \cite{CanWakBoy} that such reweighting methods find sparse solution much faster with significantly reduced errors. 

\subsection{Sparsifying multilevel transforms}

In the previous section, we discussed how the concept of sparsity is used to recover signals from possibly highly undersampled data. However, in many applications the signals are not directly sparse, but only  after the application of certain transforms. Such so-called \emph{sparsifying transforms} are often build upon systems that are equipped with a multiscale structure. Indeed, very recently a new direction of compressed sensing has been developed that is very much motivated by the sparsity structure of multiscale systems \cite{AdcHanPooRom}. Typical examples of such multiscale transforms are the \emph{wavelet} transform \cite{Dau}, the \emph{shearlet} transform \cite{GuoKutLab2006, KitKutLim, Lim}, and the \emph{curvelet} transform \cite{CanDon}. For the numerical results of this paper we will only consider the first two transforms depending on the particular signals that are to be recovered.  In particular, both, wavelet reconstructions as well as shearlet reconstructions from an incomplete amount of  Fourier measurements have been analyzed in the literature, for instance, in \cite{AdcHanKutMa} and \cite{Ma2015a, KutLimFourier}, respectively. However, the methodology of algorithm applies to any other multiscale transform other than wavelets and shearlets.


The multiscale structure of a wavelet basis comes from the use of dyadic scaling matrices of the form
\begin{align*}
    D_j = \begin{pmatrix} 2^j & 0 \\ 0 & 2^{j} \end{pmatrix}, \quad j =0,1,\ldots.
\end{align*}
Using these matrices together with  simple translations one can eventually obtain \emph{orthonormal bases} for $L^2(\R^2)$, the space of square integrable functions, of the form
\begin{align*}
\{ \phi(\cdot - m ) \, : \, m \in \Z^2 \} &\cup\{ \psi^1(D_j \cdot - m ) \, : \, j \geq 0, m \in \Z^2 \}\\
&\cup\{ \psi^2(D_j \cdot - m ) \, : \, j \geq 0, m \in \Z^2 \}\\
&\cup\{ \psi^3(D_j \cdot - m ) \, : \, j \geq 0, m \in \Z^2 \},
\end{align*}
where $\phi, \psi^1, \psi^2, \psi^3 \in L^2(\R^2)$ with certain regularity properties, see \cite{Dau} for more details.

The multiscale structure of shearlets are obtained by the use of parabolic scaling matrices of the form 
\begin{align*}
A_{2^j} = \begin{pmatrix} 2^j & 0 \\ 0 & 2^{j/2} \end{pmatrix}, \quad \widetilde{A}_{2^j} = \begin{pmatrix} 2^{j/2} & 0 \\ 0 & 2^{j} \end{pmatrix}, \quad j \in \N.
\end{align*}
In addition to the parabolic scaling matrix, a shearlet system is equipped with a directional component that can be obtained by using shear matrices
\begin{align*}
S_k = \begin{pmatrix} 1 & k \\ 0 & 1 \end{pmatrix}, \quad k \in \Z.
\end{align*}
The so-called \emph{cone-adapted shearlet system} is then defined as follows.
\begin{deff}
Let $\varphi, \psi, \psitilde \in L^2(\R^2)$. Then we call $\Phi(\phi,c) \cup \Psi(\phi, c) \cup \widetilde{\Psi}(\psitilde,c)$ a \emph{cone-adapted shearlet system}, where
\begin{align*}
\Phi(\phi,c) &= \{ \phi_m \, : \, m \in \Z^2\}, \\
\Psi(\psi,c) &= \{ \psi_{j,k,m} \, : \, j \geq 0, |k| \leq 2^{j/2}, m \in \Z^2\}, \\
\widetilde{\Psi}(\psitilde,c) &= \{ \psitilde_{j,k,m} \, : \, j \geq 0, |k| \leq 2^{j/2}, m \in \Z^2\},
\end{align*}
and
\begin{align*}
\phi_m &= \phi(\cdot - c_1m), \\ 
\psi_{j,k,m} &= 2^{3/4j} \psi(S_k A_{2^j} \cdot - cm), \\
 \psitilde_{j,k,m} &= 2^{3/4j} \psitilde(S_k^T \widetilde{A}_{2^j} \cdot - \widetilde{c}m),
\end{align*}
with $c = (c_1,c_2)^T \in \R^2_+, \widetilde{c} = (c_2,c_1)$ and the multiplication of $c$ and $\widetilde{c}$ with $m$ to be understood componentwise.
\end{deff}

\begin{figure}[H]
\includegraphics[width=0.24\textwidth]{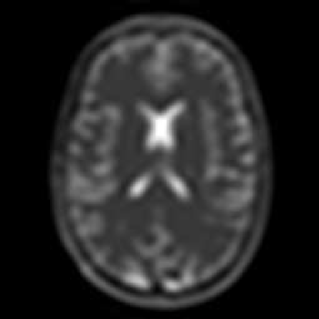}
\includegraphics[width=0.24\textwidth]{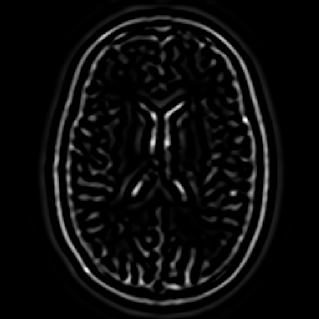}
\includegraphics[width=0.24\textwidth]{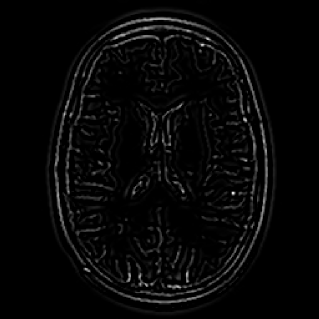}
\includegraphics[width=0.24\textwidth]{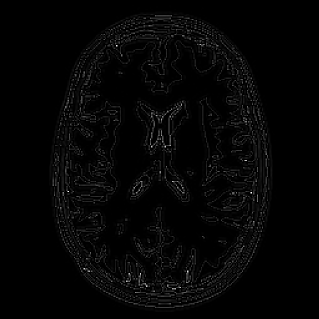}
\caption{Shearlet coefficients at different scales of the GLPU Brain phantom \cite{GLPU}. For better visual difference the contrast at scale 1,2, and 3 has been changed. 
\textbf{First}: Shearlet coefficients at scale $j = 0$. \textbf{Second}: Shearlet coefficients at scale $j =1$. \textbf{Third}: Shearlet coefficients at scale $j =2$. \textbf{Fourth}: Shearlet coefficients at scale $j =3$. } \label{Sparsity2}
\end{figure}

The shearlet coefficients $\{ (\langle u , \phi_m\rangle )_m\} \cup \{ (\langle u, \psi_{j,k,m} \rangle)_{j,k,m} \} \cup \{ (\langle u, \psitilde_{j,k,m} \rangle)_{j,k,m} \}$  computed from a shearlet transform are known to have a fast decay which ensures within the model of so-called \emph{cartoon-like images} \cite{DonSparse} an optimal sparse approximation rate \cite{KutLim}. Wavelets on the other hand, do not fulfill this optimal approximation rate. In Figure \ref{Sparsity2} we depicted some of the shearlet coefficients of each scale for an MRI test image.

\subsection{Example: Multilevel iterative reweighting for inpainting}\label{sec:inpainting}

\begin{figure*}
 \includegraphics[width=0.247\textwidth]{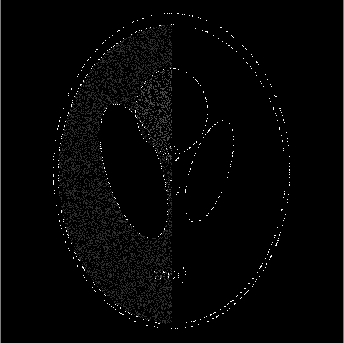}
\includegraphics[width=0.37\textwidth]{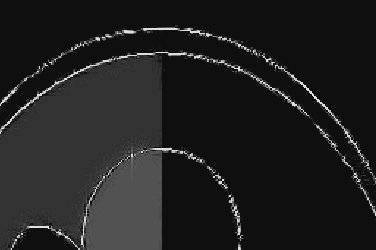}
\includegraphics[width=0.37\textwidth]{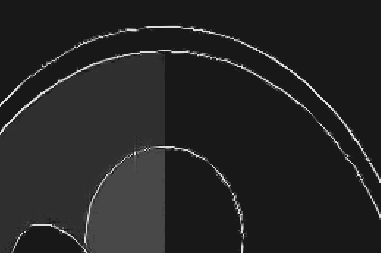}
\caption{\textbf{Left:} Sampled data of an image. \textbf{Middle:} Reconstructed image using iterative reweighting (Strategy \eqref{eq:Strat1}). \textbf{Right:} Reconstructed image using multilevel iterative reweighting (Strategy \eqref{eq:Strat2}). Note that the fine details are much better restored in the second case.} \label{fig:inpainting}
\end{figure*}

As already outlined above, reweighted $\ell^1$-minimization can significantly improve the reconstruction of certain sparse signals. The reader might wonder at this point, why an adaptive strategy for choosing the weights is beneficial if multilevel transforms are used. We demonstrate the basic idea of multilevel iterative reweighting for a simple \emph{inpainting} or \emph{image restoration} example using shearlets. Suppose $A$ denotes a masking operator, $y = Au$ is the masked version of the original image $u$, cf. Figure \ref{fig:inpainting}, and $\Psi, \Psi^{-1}$ denote the forward and backward shearlet transform, respectively.  We compute two reconstructions using the iterative hard thresholding  Algorithm \ref{algo:IHT} for two different thresholding strategies. Note that thresholding algorithms are classical and well known in the literature for such image restoration tasks. It has also been used in combination with shearlets in 
\cite{KutLimRei}.

In the following, let $\delta = (\delta(1), \ldots, \delta(N))^T \in \R^N$, where $N$ is the length of the transform coefficient vector $\Psi u$. Further, let $T_\delta: \R^N \longrightarrow \R^N$ be the hard thresholding operator applied entrywise, that is
\begin{align*}
T_\delta(c) = h, \quad \text{where} \quad h(k) = \begin{cases} c(k) & |c(k)| > \delta(k), \\ 0 & |c(k)| \leq \delta(k) \end{cases}
\end{align*}
for $k = 1, \ldots, N$.
\begin{algorithm}
\caption{Iterative hard thresholding}\label{algo:IHT}
\begin{algorithmic}
\STATE{\underline{Input}:
Measurements $y$, sampling operator $A$, sparsifying transform $\Psi$, initial values for $\delta$ and $\lambda$, thresholding strategy \eqref{eq:Strat1}  or \eqref{eq:Strat2} denoted by $f$, factor $\sigma <1$.\\
\underline{Initialization}:\\
$u_{\text{rec}},u_{\text{res}} \gets 0$;
}
\FOR{$i =1:N$}
\STATE{$u_{\text{res}} = A(y - u_{\text{rec}})$;\\
$u_{\text{rec}} = {\Psi}^{-1}(T_\delta(\Psi(u_{\text{res}} + u_{\text{rec}})))$; \\
$\delta = \sigma \cdot f(u_{\text{res}} + u_{\text{rec}}, \lambda)$;
}
\ENDFOR
\end{algorithmic}
\end{algorithm}

The shearlet coefficients $\{ (\langle u, \phi_m\rangle )_m\} \cup \{ (\langle u, \psi_{j,k,m} \rangle)_{j,k,m} \} \cup \{ (\langle u, \psitilde_{j,k,m} \rangle)_{j,k,m} \}$  can be divided into corresponding levels $j$. For the first cone this looks as follows
\begin{align*}
 (\langle u, \psi_{j,k,m} \rangle)_{j,k,m} = ((\langle u, \psi_{j,k,m} \rangle)_{(j,k,m) \in I_j})_j
\end{align*}
where $I_j$ contains all indices at scale $j$. Two different thresholding strategies involving reweighted $\ell^1$ are now for example
\begin{align}
 f_1(u,\lambda) &= \left (\frac{\lambda}{|\langle u, \psi_{j,k,m} \rangle| + \varepsilon}\right )_{j,k,m}  \label{eq:Strat1}
\end{align}
and 
\begin{align} 
 f_2(u,\mu) &= \left (\frac{\mu\max\{(|\langle u, \psi_{j,k,m}\rangle|)_{(j,k,m) \in I_j}\}}{|\langle u, \psi_{j,k,m} \rangle| + \varepsilon}\right )_{j,k,m},   \label{eq:Strat2}
\end{align}
where $\lambda, \mu$ and $\varepsilon>0$ are fixed parameters. Notice that the first strategy \eqref{eq:Strat1} is a direct application of reweighted $\ell^1$ to the thresholding based recovery algorithm Algorithm \ref{algo:IHT}, whereas strategy \eqref{eq:Strat2} involves an additional compensation factor to adapt the idea of reweighting to the multilevel structure. Such adaptations are useful as the transform coefficients of a multilevel system have a natural decrease due to scaling. Independently of the usage of reweighting, such types of compensation factors should be used in thresholding based algorithms. This is, for instance, used in \cite{HauMa} and we will discuss other choices in Section \ref{sec:Section3}.

In Figure \ref{fig:inpainting}  one can observe that by introducing level adapted weights, the information across all scales will be treated equally important. Without such weights, the fine detail coefficients are underrated in terms of their importance which might lead to crucial quality loss of detail information.

Algorithm \ref{algo:IHT} is simple and effective, however, it is not very sophisticated and also not sufficient for our purposes as in this work we also want to consider an additional regularizer. Hence, we continue with an introduction to \emph{split Bregman}.

\subsection{Split Bregman}
\label{SB}
Split Bregman is an effective algorithm to solve constrained optimization problems by introducing a split variable and solving the resulting decoupled problem with \emph{Bregman Iterations}. It was proposed in 2009 by Goldstein and Osher and became a popular method since then \cite{SB,SB2,Wu,Plonka,Steidl}. However, as we have already mentioned above, the algorithm was already invented under the name ADMM in \cite{GabMer,EckBer}; see also \cite{Setzer,Esser}. 

Even though split Bregman can handle general convex regularizers it is sufficient for our purposes to focus on $\ell^1$-regularized problems in the analysis formulation. We briefly follow the steps in \cite{Burger,SB,Yin} to derive the basic form of the algorithm.

Consider the basis pursuit problem in analysis formulation 
\begin{align}
    \label{BPinAna}
    \min_u  \norm{W \Psi u}_1 \quad \mbox{ subject to } \quad  Au = y,
\end{align}
for a possibly redundant dictionary $\Psi \in \R^{N}\times n, N, n \in \N$, a measurement matrix $A \in \C^{m\times n}$, and a diagonal weighting matrix with entries
$
W(l,l)$ for $l=1,\dots,N$.
Then instead of using a continuation method for enforcing the constraint, i.e., taking $\beta \to \infty$ in
\begin{align*}
    u= \argmin_u \norm{W\Psi u}_1 + \frac{\beta}{2} \norm{y -Au}_2^2,    
\end{align*}
problem \eqref{BPinAna} is transformed into a sequence of unconstrained problems using \emph{Bregman iterations}
\begin{align}
\label{BPinAnaBreg}
\begin{cases}
    u_{k+1} = \argmin \limits_u \norm{W \Psi u}_1 + \frac{\beta}{2} \norm{Au - y_k}_2^2,\\
    y_{k+1} = y_k + y - Au_{k+1},
\end{cases}
\end{align}
for a suitable $\beta >0$. For keeping the presentation in this theoretical study as concise as possible we focus  on  noiseless measurements and remark that if noise is present, model \eqref{BPinAna} and the  iterations \eqref{BPinAnaBreg} have to be adapted carefully, cf. \cite{Yin, Burger}. A performance of the presented framework in this situation is demonstrated in an upcoming work \cite{MaMaeFunSchKutSchKol}. 
We continue by introducing a split variable $d = \Psi u$ for the $\ell^1$-part of the minimization problem in \eqref{BPinAnaBreg} and executing an additional Bregman iteration step to obtain
\begin{alignat*}{1}
\begin{cases}
(u_{k+1},d_{k+1}) & = \argmin \limits_{u,d}  \norm{Wd}_1   + \frac{\beta}{2} \norm{Au - y_k}_2^2  \\
& \qquad  \qquad + \frac{\mu}{2} \norm{d - \Psi u - b_k}_2^2, \\ 
b_{k+1} & =   b_k + \Psi u_{k+1} - d_{k+1},\\
y_{k+1} & =  y_k + y - Au_{k+1}.
\end{cases}
\end{alignat*}
To solve the $(u,d)$-minimization problem one or multiple \emph{nonlinear block Gauss-Seidel} iterations are used, which alternate between minimizing with respect to $u$ and $d$. This yields the \emph{split Bregman Algorithm}
\begin{algorithm}
\caption{Split Bregman algorithm}\label{SBAlgo}
\begin{algorithmic}

\FOR{$i =1:N$}
\STATE{$u_{k+1} =\argmin_{u}  \frac{\beta}{2} \norm{Au - y_k}_2^2 + \frac{\mu}{2} \norm{d^{\text{cur}} - \Psi u - b_k}_2^2$\\
$d_{k+1}  = \argmin_d \norm{Wd}_1 + \frac{\mu}{2} \norm{d - \Psi u^{\text{cur}} - b_k}_2^2,$
}
\ENDFOR
\STATE{$b_{k+1}  =   b_k + \Psi u_{k+1} - d_{k+1}$,\\
$y_{k+1}  =  y_k + y - Au_{k+1},$}
\end{algorithmic}
\end{algorithm}

where $x^{\text{cur}}$ denotes the latest available stage of the variable $x \in \{ d, u \}$. Note that the solution of the $d$-subproblem  is explicitly given by \emph{soft-thresholding}
\[
d_{k+1}(l) = \text{shrink}\left(\left(\Psi u^{\text{cur}}\right)(l) + b_{k}(l), \frac{W(l,l)}{\mu} \right), 
\]
for $l=1,\dots,N$ where
\[
\text{shrink}\left(z,\lambda\right) = \begin{cases} \max\left(\norm{z} - \lambda,0\right) \frac{z}{\norm{z}}, & z \neq 0, \\
0, & z = 0. \end{cases}
\]
In \cite{SB} it was furthermore observed that the minimization with respect to $u$ in \eqref{SBAlgo} does not have to be solved to full precision and in many applications only few steps of an iterative method are sufficient. 

\subsection{TV and TGV}\label{TGV}

\emph{Total Variation (TV)} based methods were initially proposed by Rudin, Osher, and Fatemi in 1992  for image denoising \cite{RudOshFat} and are now widely used for image reconstruction and compressed sensing, see for example \cite{RecPF,NeedellWard}. TV is based on the assumption that the reconstructed image is piecewise constant and therefore gradient sparse. This results in preserving sharp edges. But for realistic images, which are usually not piecewise constant, this can lead to severe \emph{oil painting artifacts} or  \emph{staircasing effects} leading to unnatural looking reconstructed images. 

\emph{Total Generalized Variation} (TGV) is a generalization of TV and has been proposed to improve on these issues by involving higher order derivatives \cite{BreKunPoc}. We will now briefly give a definition of the second order TGV regularizer in $\R^2$ together with some basic facts. Its general derivation and more details on this subject can be found in \cite{TGVinverse,Bredies2015,BreKunPoc,Holler,Holler2}. Please note that the discretized TGV-model, which is solved numerically, was already introduced in \cite{SS08}.

The so-called \emph{pre-dual} formulation of second order TGV is given by
\begin{align}
\TGV_{\alpha}^2(u) =  \sup  \bigg\{ & \int_{\Omega}  u \divergenz^2 v \, dx : v \in C_c^2(\Omega,\mathcal{S}^{2\times 2}), \nonumber \\
& \norm{v}_{\infty} \leq \alpha_0, \norm{\divergenz v}_{\infty} \leq \alpha_1 \bigg\},
\end{align}
for $\alpha = (\alpha_0,\alpha_1) \in \R_+^2$, $\Omega \subseteq \R^2$ a bounded Lipschitz domain, $\mathcal{S}^{2 \times 2}$ the space of symmetric $2\times 2$ matrices, and $u \in L^1(\Omega,\C)$. Thereby the divergences are defined as
\[
(\divergenz v)_i =\sum_{j=1}^2  \frac{\partial v_{ij}}{\partial x_j}, \qquad i=1,2, 
\]
and 
\[
\divergenz^2 v = \sum_{i,j=1}^2 \frac{\partial^2 w_{ij}}{\partial x_i \partial x_j},
\]
together with the norms
\[
\norm{v}_{\infty} = \sup_{l\in\Omega} \left(\sum_{i,j=1}^2 |v_{ij}(l)|^2 \right)^{1/2},
\]
and
\[
\norm{\divergenz v}_{\infty} = \sup_{l\in\Omega} \left(\sum_{i=1}^2 |(\divergenz v)_i(l)|^2 \right)^{1/2}. 
\]
Under certain conditions, an equivalent and more convenient form of $\TGV_{\alpha}^2$ is given by the \emph{minimum representation} as
\begin{align}
\TGV_{\alpha}^2(u) = \inf_{v \in \text{BD}(\Omega,\C^2)} \alpha_1 \norm{\nabla u -v}_1  + \alpha_0 \norm{\mathcal{E}(v)}_1,
\end{align}
where $\text{BD}(\Omega, \C^2)$ is the space of symmetric tensor fields of bounded deformation and $\mathcal{E}$ the symmetrized derivative defined as
\[
\mathcal{E}(v) = \left(\begin{matrix}\partial_x v_1 &  \frac{1}{2}(\partial_y v_1 +  \partial_x v_2) \\ \frac{1}{2}(\partial_y v_1 + \partial_x v_2) & \partial_y v_2  \end{matrix}\right).
\]
In this form $\TGV_{\alpha}^2$ can be interpreted as balancing the first and second derivatives of $u$ controlled by the ratio of $\alpha_0$ and $\alpha_1$. In \cite{BreKunPoc} and \cite{KnoBrePocRud} it was observed that the use of TGV as a regularizer indeed leads to reconstructed images with sharp edges but without the staircaising effects of TV. Being furthermore convex and lower semi-continuous with respect to $L^1$-convergence, this makes TGV a numerical feasible and therefore a suitable alternative for TV \cite{TGVinverse}.

\section{Proposed multilevel based reweighting algorithm with TGV}\label{sec:Section3}

\subsection{Model and discretization}
In this section we aimto develop an algorithm for solving the multilevel reweighted $\ell^1$-problem. In order to do so, the split Bregman approach introduced in Section \ref{SB} shall be equipped with an appropriate iteratively reweighted soft-thresholding procedure.  For a further reduction of artifacts and improving the reconstruction of piecewise constant as well as smooth regions we not only use the regularizer belonging to the reweighted multilevel decomposition but also TGV as a second regularization term. 

Let $A$ be a measurement operator, let $y$ the measurements of our signal of interest $u$, and let $\sigma>0$ be a fidelity parameter. The recovery problem can then be stated as 
\begin{align*}
    \min_u \sum_{j=1}^{\infty} \lambda_j \norm{W_j \Psi_j u}_1 + \TGV_{\alpha}^2(u)  \quad \mbox{ subject to } \quad Au = y, 
\end{align*}
where $\Psi_j$ corresponds to the $j$-th subband of the multilevel transform $\Psi $, $\lambda_j$ are regularization parameters accounting for the multilevel structure of $\Psi $ and $W_j$ are diagonal matrices containing the weights. For the sake of clearness we have assumed that there is only one subband per level otherwise an additional index has to be attached to $\Psi_j$ to specify the current subband.
Note that after we have established a basic split Bregman framework for solving the minimization problem we will aim to update $\lambda_j$ and $W_j$ iteratively.
Using the characterization of $\TGV_{\alpha}^2$ presented in Section \ref{TGV} the objective can be rewritten as 
\begin{align}
\label{Objective}
\min_{u,v} \sum_{j=1}^{\infty} \lambda_j \norm{W_j \Psi_j u}_1 + \alpha_1 \norm{\nabla u - v}_1 + \alpha_0 \norm{\mathcal{E}(v)}_1.
\end{align}
For the \emph{discretization}, let $u\in\C^{n^2}$ be the vectorized finite-dimensional image of interest which is for simplicity assumed to be of square size. Let $A\in \C^{m\times n^2}$ be the finite dimensional measurement matrix and $y\in \C^m$ the observed data. Furthermore, let $\nabla^f$ and $\nabla^b$ denote a discrete gradient operator with periodic boundary conditions using forward and respectively backward differences. Following \cite{BreKunPoc,Holler2} we  approximate the derivatives in \eqref{Objective} by
\[
\nabla u \approx \nabla^f u = \left(\begin{matrix} \nabla_x^f u\\ \nabla_y^f u \end{matrix}\right)
\]
and 
\[
\mathcal{E}(v) \approx \mathcal{E}^b v =  \left(\begin{matrix} \nabla_x^b v_x & \frac{1}{2}(\nabla_y^b v_x + \nabla_x^b v_y )\\ \frac{1}{2}(\nabla_y^b v_x + \nabla_x^b v_y ) & \nabla_y^b v_y \end{matrix} \right).
\]
A finite dimensional approximation of \eqref{Objective} is then given by
\begin{align}
\label{discrObj}
\min_{u,v} \sum_{j=1}^{J} \lambda_j \norm{W_j \Psi_j u}_1 + \alpha_1 \norm{\nabla^f u - v}_1 + \alpha_0 \norm{\mathcal{E}^bv}_1,
\end{align}
where $J$ is some fixed a priori chosen maximum scale and $\Psi$ is the discrete transform acting on the vectors. 
Let us furthermore introduce the notation 
\begin{align}
\begin{split}
	\Psi u  = \left(\Psi_j  u \right)_{j=0,\dots,J}  = (\langle \psi_{j, l}, u \rangle)_{j=0,\ldots, J, l =1, \ldots, N_j} \label{eq:CoeffNj}
	\end{split}
\end{align}
for dividing the analysis coefficients into $J$ subbands, each consisting of $N_j \in \N$ elements.
For a documentation of the discrete transforms, we refer the interested reader to Chapter 8 in \cite{Walnut} for wavelets and \cite{KutLimRei} for shearlets. 
Note that the $\ell^1$-norm in the second summand of \eqref{discrObj} is thereby defined as 
\[
\norm{v}_1 = \sum_{l=1}^{n^2} \left(|v_x(l)|^2 +|v_y(l)|^2 \right)^{1/2},
\]
and for the third summand as
\[
\norm{e}_1 = \sum_{l=1}^{n^2} \norm{e(l)}_F = \sum_{l=1}^{n^2} \norm{\left(\begin{matrix}
 e(l)_1 & e(l)_2 \\
 e(l)_2 & e(l)_3
\end{matrix} \right)}_F, 
\]
where $\norm{\cdot}_F$ is the Frobenius norm of a $2 \times 2$ matrix.
\subsection{Split Bregman framework}
The proposed constrained optimization problem can be casted into the form given in  \eqref{BPinAna} by introducing the variable ${{\bf u} = (u,v)^T}$, together with the matrix
\[
{\bf \Psi } = \left(\begin{matrix} \Psi  & 0 \\  \nabla^f & - I \\ 0 &  \mathcal{E}^b  \end{matrix} \right).
\] 
In order to  come up with the explicit form of the resulting split Bregman algorithm as given in Section \ref{SB}, let us  split  as follows:
\[
\left( \begin{matrix}
w\\d \\t 
\end{matrix}
\right) 
= \left(
\begin{matrix}
 \Psi u \\ \nabla^f  u -v  \\ \mathcal{E}^bv
\end{matrix}  \right). 
\]
The $(u,v)$-subproblem of Algorithm \ref{SBAlgo} is then given by 
\begin{align}
\begin{split}
\label{uProb}
 &(u_{k+1},v_{k+1}) \\
 &=  \argmin_{u,v}  \frac{\beta}{2} \norm{Au - y_k}_2^2   + \frac{\mu_1}{2} \norm{w^{\cur}- \Psi u - b^w_k}_2^2 \\ 
   &+ \frac{\mu_2}{2}\norm{d^{\cur} - (\nabla^f u - v) - b^d_k }_2^2   + \frac{\mu_3}{2}\norm{t^{\cur} - \mathcal{E}^b v - b^t_k}_2^2.
\end{split}
\end{align}
We furthermore obtain the subproblems
\begin{align}
\label{subw}
w^j_{k+1} = \argmin_{ w^j} \lambda_j  \norm{W_j  w^j}_1 
+ \frac{\mu_1}{2}  \norm{w^j - \Psi_j u^{\cur} - b^{w,j}_k}_2^2,
\end{align}
for each subband $j=1,\dots,J$, as well as
\begin{align}
\label{subd}
d_{k+1} = \argmin_d \alpha_1 \norm{d}_1 + \frac{\mu_2}{2} \norm{d - (\nabla^f u^{\cur} - v^{\cur}) - b^{d}_k}_2^2,
\end{align}
and 
\begin{align}
\label{subt}
t_{k+1} = \argmin_t \alpha_0 \norm{t}_1 + \frac{\mu_3}{2} \norm{t - \mathcal{E}^b v^{\cur} - b^t_k}_2^2.
\end{align}
Note that the regularization parameters $\lambda_j, \alpha_0,$ and $\alpha_1$ have thereby been subsumed into a weighting matrix $W$. Also we are allowing some more flexibility by incorporating different values for $\mu_i$ for $i=1,2,3$. 
Furthermore we obtain the following \emph{Bregman updates:}
\begin{align*}
\begin{cases}
b^w_{k+1} & = b^w_k + \Psi u_{k+1} - w_{k+1}. \\ 
b^d_{k+1} & = b^d_k + (\nabla^f u_{k+1} -v_{k+1}) -d_{k+1},\\
b^t_{k+1} & = b^t_k + \mathcal{E}^b v_{k+1} - t_{k+1},
\end{cases}
\end{align*}
as well as
\begin{align*}
    y_{k+1}  = y_k + y - Au_{k+1}.
\end{align*}

\subsection{Solutions of the subproblems}
The solution of the subproblem \eqref{uProb} can be obtained by setting the first derivatives with respect to $u$, $v_x$, and $v_y$ to zero. We then obtain the linear system
\begin{align}
\label{sys}
   \left(\begin{matrix}
    b_1 & b_4^* & b_5^* \\
    b_4 & b_2 & b_6^* \\
    b_5 & b_6 & b_3
   \end{matrix} \right) 
   \left(
   \begin{matrix}
    u \\ v_x \\ v_y
   \end{matrix}
   \right)
   = \left(
   \begin{matrix}
    R_1 \\ R_2 \\ R_3
   \end{matrix}
   \right),
\end{align}
where $b_i$ are $n^2\times n^2$ block matrices defined as
\begin{align*}
     &  b_1 = \beta A^*A + \mu_1 \Psi^* \Psi + \mu_2 (\nabla^f)^* \nabla^f, \\
    & b_2 = \mu_3 (\nabla_x^b)^* \nabla_x^b + \frac{\mu_3}{2}   (\nabla_y^b)^* \nabla_y^b + \mu_2 I,\\
    & b_3 = \mu_3 (\nabla_y^b)^* \nabla_y^b + \frac{\mu_3}{2} (\nabla_x^b)^* \nabla_x^b
     + \mu_2 I, \\
     & b_4 = - \mu_2  \nabla_x^f, \\
     & b_5 = -\mu_2 \nabla_y^f,\\
     & b_6 = \frac{\mu_3}{2} (\nabla_x^b)^* \nabla_y^b,
\end{align*}
and the components of the right hand side are given by
\begin{alignat*}{2}
     R_1 &= \beta A^*(y + y_k)  + \mu_1 \Psi^*(w^{\cur} - b_k^w) \\
     &\qquad + \mu_2 (\nabla^f)^* (d^{\cur} -b_k^d), \\
     R_2 &=   \mu_2 ( b_{k,x}^d - d_x^{\cur}) + \mu_3 \big[  (\nabla_x^b)^* (t_1^{\cur}-b^t_{k,1})\\
          &\qquad +(\nabla_y^b)^* (t_2^{\cur}-b^t_{k,2})\big], \\
  R_3 &=   \mu_2 ( b_{k,y}^d - d_y^{\cur})   + \mu_3 \big[  (\nabla_x^b)^* (t_2^{\cur}-b^t_{k,2})\\
 & \qquad+(\nabla_y^b)^* (t_3^{\cur}-b^t_{k,3})\big].
\end{alignat*}
Similar to \cite{SB}, it was observed in \cite{GuoQinYin}, that in many cases the linear system in \eqref{sys} can be efficiently solved by using the $2$D-Fourier transform $\mathcal{F} \in \C^{n^2 \times n^2}$. Note that $\nabla^f$ and $\nabla^b$ are circulant, since they correspond to periodic boundary conditions. Therefore 
\[
\mathcal{F} \nabla^* \nabla \mathcal{F}^*
\]
is a diagonal matrix. For a tight frame $\Psi$ we have 
\begin{align}
\Psi^* \Psi = a I,\label{eq:TightFrame}
\end{align}
where $a\in\R$ is the frame bound of $\Psi$. This is for example the case in \cite{GuoQinYin}, where the Fast Finite Shearlet Transform (\cite{Hae13,Hae14}) was used, which forms a Parseval frame for $\R^{n\times n}$, i.e., Equation \eqref{eq:TightFrame} with $a = 1$. Note that also for the non-tight shearlet system of ShearLab \cite{KutLimRei}  the matrix $\mathcal{F}\Psi^*\Psi\mathcal{F}^*$ is diagonal and can be explicitly computed.

In the case of  subsampled Fourier measurements, the measurement matrix can be written as
\begin{align*}
A = P \mathcal{F},
\end{align*}
where $P \in \left\{0,1\right\}^{m\times n^2}$ is selecting or discarding measurements. In this case 
\[
A^* A = \mathcal{F}^* P \mathcal{F}
\]
is naturally diagonalized by the 2D Fourier transform. 

If all blocks $b_i$ for $i=1,\dots,6$ can be diagonalized in this way, the authors of \cite{GuoQinYin} proposed to multiply with a preconditioner matrix from the left to obtain the system
\begin{align}
\label{diagsys}
\left( \begin{matrix}
\widehat{b_1} & \widehat{b_4}^* &\widehat{b_5}^*\\
\widehat{b_4} & \widehat{b_2} & \widehat{b_6}^* \\
\widehat{b_5} & \widehat{b_6} & \widehat{b_3} \end{matrix} \right)
\left( 
\begin{matrix}
\mathcal{F}u \\ \mathcal{F} v_x \\ \mathcal{F} v_y \end{matrix} \right) 
 = \left( \begin{matrix}\mathcal{F}R_1 \\ \mathcal{F}R_2 \\ \mathcal{F} R_3 \end{matrix} \right),
\end{align}
where each $\widehat{b_j} = \mathcal{F} b_j\mathcal{F}^*$ is a $n^2\times n^2$ diagonal matrix. A closed form solution can then be obtained by applying Cramer's rule.   

In other applications, such as inpainting or reconstruction from the subsampled cosine transform, the measurement matrix can be written as
\begin{align}
A = P \mathcal{T}, \label{eq:AUnit}
\end{align}
for a unitary matrix $\mathcal{T} \in \C^{n^2 \times n^2}$. In the ADMM model of \cite{GuoQinYin} it is proposed to include an additional split to deal with the fact that $A^* A$ cannot be diagonalized under $\mathcal{F}$. 

\noindent However, note that a representation of the form \eqref{eq:AUnit} is not possible for instance in 
\begin{enumerate}
\item \emph{Partial Parallel Imaging (PPI)} \cite{PPI,PPI2} in MRI,  where the image is to be reconstructed by using subsampled Fourier data from parallel scans of multiple coils. 
\item \emph{Computed tomography}, where the image has to be reconstructed from \emph{fan beam} projections, see Section \ref{sec:Radon} for more details.
\end{enumerate}
In such cases we therefore propose to solve the linear system of equation \eqref{sys} only approximately. As it was observed in \cite{SB}, usually only a few steps of an iterative solver are sufficient for the convergence of the resulting algorithm. A proof of this fact can be found in \cite{Nien}. In order to save memory and computation time, we still want to use  that all $d_j$ are diagonal under $\mathcal{F}$ except for $\widehat{d_1}$.  Therefore we multiply with the same preconditioner matrix as above and obtain a similar system as in equation \eqref{diagsys}, with the difference that $\widehat{d_1}$ is not diagonal anymore. With a few explicitly given steps the block matrix can be brought to \emph{lower triangular form}. In this way 
we only need to apply an iterative solver to one $n^2 \times n^2$ system involving $A^* A$ instead of solving the entire block system.  
In both of the above listed cases, we are using a projected conjugate gradient method together with a warm start (obtained through the previous iterations of the  split Bregman algorithm), so that only very few steps are necessary for sufficient precision; see also Section \ref{sec:Radon} and \cite{MaMaeFunSchKutSchKol}. 

Let us now briefly discuss the solution of the other subproblems: 
As described in Section \ref{SB} a closed-form solution of \eqref{subw} is given by
\begin{align}
\label{shrinkw}
w_{k+1}^j(l) = \text{shrink}\left(\left(\Psi_j u^{\cur}\right)(l) + b_k^{w,j}(l), \frac{\lambda_j  {W_j}(l)}{\mu_1} \right), 
\end{align}
for $l=1,\dots,N_j-N_{j-1}+1$. 
For equation \eqref{subd} we obtain
\begin{align*}
d_{k+1}(l) = \text{shrink}_2\left( \nabla^f u^{\cur}(l) - v^{\cur}(l) + b_k^d(l),\frac{\alpha_1}{\mu_2} \right),
\end{align*}
for $l=1,\dots,n^2$ and the shrinkage rule
\begin{align*}
\text{shrink}_2\left(x,\lambda\right) = \begin{cases}
\max(\norm{(x}_2 - \lambda),0) \frac{x}{\norm{x}_2}, & x \neq 0,\\
0, & x= 0.
\end{cases}
\end{align*}
Similarly, the solution of \eqref{subt} is given by
\begin{align*}
t_{k+1}(l) = \text{shrink}_F\left(\left(\mathcal{E}^b v^{\cur}\right)(l) + b^t_k(l), \frac{\alpha_0}{\mu_3} \right),
\end{align*}

for $l=1,\dots,n^2$ and 
\begin{align*}
\text{shrink}_F\left(x,\lambda\right) = \begin{cases}
\max(\norm{(x}_F - \lambda),0) \frac{x}{\norm{x}_F}, & x\neq 0,\\
0, & x = 0.
\end{cases}
\end{align*}

\subsection{Combining reweighted $\ell^1$ with multiscale transforms}\label{Sec:IIIb)}

As we already outlined in Section \ref{sec:Section2}, the idea of \emph{reweighted $\ell^1$} as proposed by Cand\`{e}s et al. in \cite{CanWakBoy} is to improve the $\ell^1$-norm as a sparsity regularizer by an iterative reweighting of the nonzero coefficients. In contrast to $\ell^1$-minimization, $\ell^0$ only counts the number of nonzero coefficients and thus spares larger coefficients more than $\ell^1$ does. 
Hence the idea of reweighted $\ell^1$ is that since larger coefficients of an iterative solution are more likely to be nonzero in the true signal, they should be multiplied with a smaller weight during the optimization process.
Put differently, the guiding principle of reweighted $\ell^1$ is that small coefficients of an iterative solution are likely going to be zero in the true signal. However, this principle is not valid for multiscale sparse signals, i.e., signals that can be sparsely represented using a multiscale transform. The magnitudes of multiscale coefficients naturally decrease with increasing scales, but the high scale nonzero coefficients of an iterative solution are not necessarily less important or more likely zero in the actual signal, if compared to low scale coefficients which are intrinsically larger. In the following section we are aiming to compensate for this misfit by including additional weighting parameters for each level in the transformation.

Suppose $u \in \R^{n^2}$ is the true signal. Let us recall the notation
\begin{align*}
\begin{split}
	\Psi u  = \left(\Psi_j  u \right)_{j=0,\dots,J}  = (\langle \psi_{j, l}, u \rangle)_{j=0,\ldots, J, l =1, \ldots, N_j}, 
	\end{split}
\end{align*}
for dividing the analysis coefficients into $J$ subbands, each consisting of $N_j$ elements. In \cite{AhmSch} a multi dictionary reweighting algorithm (Co-L1) was proposed  which iteratively updates $\lambda_j^{k}$ in
\begin{align}
\label{prob}
u_{k+1} = \argmin_{u} \sum_{j=0}^{J} \lambda_j^{k} \norm{W_j \Psi_j u}_1 \quad \text{ subject to } \quad Au = y,
\end{align}
by setting
\begin{align}
\label{Schniter}
\lambda_j^{k} = \frac{N_j}{\varepsilon + \norm{\Psi_j u^{k}}_1},
\end{align}
and $W_j = I$ for all iterations of solving \eqref{prob}. It was shown therein that the resulting algorithm can be interpreted as applying a Majorization-Minimization algorithm to the unconstrained formulation of \eqref{prob} with the regularizer
\begin{align*}
    \sum_{j=0}^J N_j \log\left(\varepsilon + \norm{\Psi_j u}_1 \right).
\end{align*}
This update rule was proposed in \cite{AhmSch} for a composition of multiple different dictionaries instead of just one multiscale dictionary divided into its subbands. In the latter case it is less likely to expect that $\log \left(\varepsilon + \norm{\Psi_j u}_1 \right)$ promotes the sparsity structure of $u$ within each of the subbands $\Psi_j$ sufficiently. Indeed, as it was argued in \cite{CanWakBoy}, the log-sum penalty is more sparsity enforcing than the $\ell^1$-norm by putting a larger penalty on small nonzero coefficients. In the case of the $\ell^1$-norm of an entire subband this approach seems to be less effective in promoting the sparsity within each level. It was furthermore proposed in \cite{AhmSch} (Co-IRW-L1, Algorithm 4) to combine the update rule \eqref{Schniter} with the classical elementwise reweighting update 
\begin{align}
\label{IRL1}
W_{j} = \diag \left( \frac{1}{\varepsilon + \left| \langle \psi_{j,l},u^{k} \rangle \right| } \right), 
\end{align}
for $j=1,\dots,J$.
However, note that this combination is very different to what we are aiming for, since there is even more emphasize put on penalizing the smaller coefficients in higher levels which can happen to delete too many highscale coefficients.

\begin{figure*}
\includegraphics[width=0.49\textwidth]{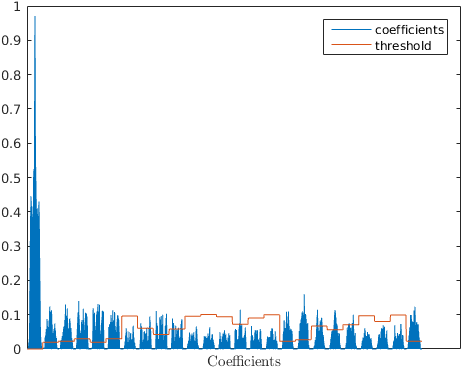}
\includegraphics[width=0.49\textwidth]{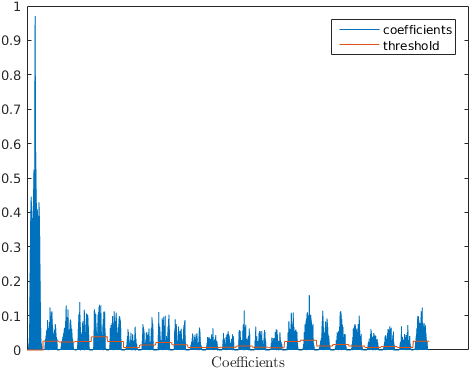}
\\[1ex]
\includegraphics[width=0.49\textwidth]{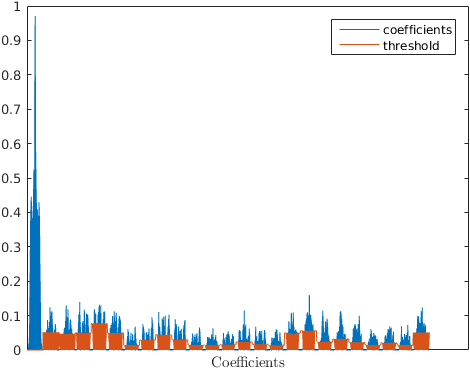}
\includegraphics[width=0.49\textwidth]{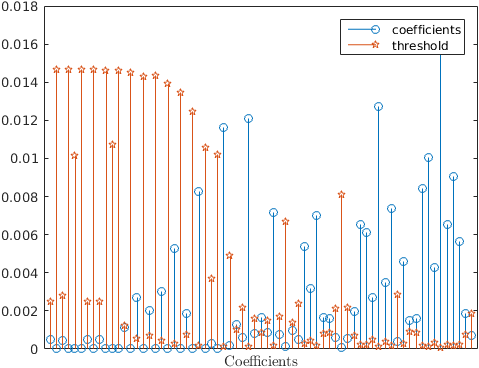}
\caption{\textbf{Upper left}: Coefficients (in blue) with a multi-level reweighting curve as in \eqref{Schniter} (in orange). \textbf{Upper right}: Coefficients (in blue) with multi-level curve (in orange) using \eqref{our} only. \textbf{Lower left}: Coefficients with the proposed reweighting strategy (in orange). \textbf{Lower right}: Zoom of proposed multilevel reweighting strategy.} \label{Reweighting}
\end{figure*}

This fact is visualized from a different point of view in Figure \ref{Reweighting}, where we have depicted the shearlet coefficients of a MRI phantom introduced in \cite{GLPU} together with the reweighting rule we have just discussed in the top-left of the figure.  The shearlet coefficients are depicted in blue and the values of $\frac{\lambda_j}{\mu_1}$ for a realistic value $\mu_1$ are shown in orange. Note that according to the update rule  \eqref{shrinkw} of the split Bregman algorithm everything below the orange curve would be thresholded. 

\subsection{Multilevel adapted reweighting}{Sec:IIIc)}
Considering the previous discussion one of the disadvantages is that the weights corresponding to higher levels might become too large. This can be prevented, for instance, by choosing the regularization parameters as 
\begin{align}
\label{our}
\lambda_j =  \max\left\{\left| \langle \psi_{j,l}, u \rangle \right| : l = 1,\dots,N_j \right\},
\end{align}
for $j =1, \ldots, J$ and zero otherwise, i.e., if $j =0$. Notice that the necessity of such weights was already elaborated in Section \ref{sec:inpainting}. 
Note that we set $\lambda_0=0$, since for real life signals the low frequency part is usually not sparse, see also Figure \ref{Sparsity2}. This was also observed in  \cite{SelFig}, where it was shown that this idea can be accomplished more effectively using an analysis instead of a synthesis prior. The attentive reader might wonder that the magnitude of the analysis coefficients could be very irregular per level, in particular, one could have strong outliers. For such cases it is better to take a quantile instead of the maximum over all coefficients. Note that this choice of $\lambda_j$ is a heuristic substitute accounting for the unknown constant in the theoretical decay of the multilevel coefficients.  A schematic representation from the thresholding perspective of split Bregman can be found in the second image of Figure \ref{Reweighting}. 

\begin{figure*}
\includegraphics[width=0.495\textwidth]{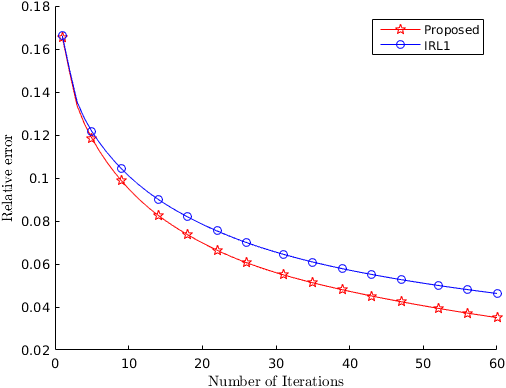}
\includegraphics[width=0.495\textwidth]{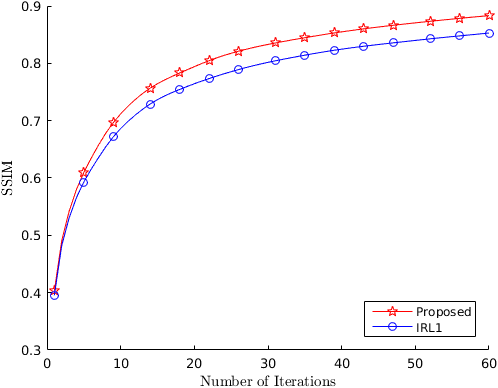}
\caption{Convergence plot: Behavior of the error and structured similarity index for reweighting with true shearlet coefficients with respect to increasing number of iterations. Used signal: phantom from \cite{GLPU}. Reconstructed from 6\% of Fourier data with proposed algorithm without TGV (see section \ref{sec:Section4} for details). In IRL1 we are choosing $W_j$ as in \eqref{IRL1} and $\lambda_j$ constant and for the proposed emthod we additionally define $\lambda_j$ as in \eqref{our}. \textbf{Left figure}: Relative error of in each iteration.  \textbf{Right figure}: SSIM of each iteration.} \label{true}
\end{figure*}

Our proposed method combines the classical reweighting of \eqref{IRL1} with the above choice for $\lambda_j$. The idea behind this is that we are still using the power of pointwise iterative reweigthing, but since our multiscale coefficients naturally come in a decreasing order in magnitude, we apply it to each level separately weighted with $\lambda_j$. That means that within each level we follow the democratic philosophy of reweighting which is that small analysis coefficients of the current iterate $\Psi_j u^{\cur}$ are likely to be zero in $\Psi_j u$. By multiplying with the latter choice of $\lambda_j$ we  also gain more control and account for the multilevel structure of $\Psi u$. An artificial experiment using the (in reality unknown) true analysis coefficients demonstrates that this update rule does seem to perform better than standard reweighting without such a compensation of multilevel weights, see Figure \ref{true}. For the signal $u$ we are choosing the phantom of \cite{GLPU}. Using the 
usually unknown shearlet coefficients for the construction of $\lambda_j$ and $W_j$ for $j=1,\dots,J$ as explained above we are reconstructing $u$ from only $6 \%$ of its Fourier measurements obtained by radial lines through the k-space origin by using the proposed algorithm of the last section, but without the additional TGV regularizer. Note that both reconstructions approximately start with the same error, which indicates that the set of regularization parameters is chosen equally good. However, it can be said that tuning the iterative reweighted shrinking method within the split Bregman framework without the automatic choice of the subband-weights $\lambda_j$ is rather difficult and highly signal dependent. This example further shows the potential of reweighting as the data is significantly undersampled.

 \subsection{Proposed algorithm}
 
Having explained the idea of our method, we now state the final resulting algorithm that is a composition of the split Bregman framework for solving the constrained optimization problem \eqref{discrObj} and the previously explained idea of multilevel weighting and iteratively reweighting, respectively.
In contrast to the traditional reweighted $\ell^1 $ approaches we do not iterate between solving the $\ell^1$-problem up to convergence and updating the weights. We propose to incorporate the multilevel adapted reweighting rule directly into the split Bregman algorithm. This is done in such a way that only the shrinking of the $w$-subproblem is changed to a \emph{multilevel adapted, iteratively reweighted shrinking rule}. Note that by choosing the level-weights $\lambda_j$ depending on the magnitude of the signal coefficients the resulting method appears to stable towards the alternation of signals. 

\begin{algorithm*}
\caption{Proposed algorithm}
\label{alg:buildtree}
\begin{algorithmic}
\STATE{\underline{Input}: \\
Measurement operator $A$, multilevel transform $\Psi$,\\
regularization parameters: $\alpha_0, \alpha_1,\mu_1,\mu_2,\mu_3,\beta,$\\
iteration numbers $N$ and maxIter.\\
\underline{Data}:\\
Measured data $y$.\\
\underline{Initialization}:\\
$ k \gets 0$;\\
$ u_0 \gets A^*y$;\\
$ y_0,v_0, d_0, w_0,t_0,b_0^t,b_0^d,b_0^w \gets 0$;
}
\WHILE{$k \leq \maxIter$}
\FOR{$i=1,\ldots,N$}
\STATE{$(u_{k+1},v_{k+1}) \gets$ \text{solve linear system \eqref{diagsys}};}
\FOR{$j=1,\ldots,J$}
\STATE{$\lambda_j =  \max\left\{\left| \langle \psi_{j,l}, u \rangle \right| : l = 1,\dots,N_j  \right\}$;\\
$W_{j} = \diag \left( \frac{1}{\varepsilon + \left| \langle \psi_{j,l},u^{\cur} \rangle \right| } \right)$;\\
$w_{k+1}^j(l) \gets \text{shrink}\left(\left(\Psi_j u^{\cur}\right)(l) +  b_k^{w,j}(l), \frac{\lambda_j W_j(l)}{\mu_1} \right), \quad l=1,\dots,N_j$};\\
\ENDFOR
\STATE{$d_{k+1}(l) \gets \text{shrink}_2\left( \nabla^f u^{\cur}(l) - v^{\cur}(l) + b_k^d(l),\frac{\alpha_1}{\mu_2} \right),\quad l=1,\dots,N_j$;\\
${t_{k+1}(l) \gets \text{shrink}_F\left(\left(\mathcal{E}^b v^{\cur}\right)(l) + b^t_k(l), \frac{\alpha_0}{\mu_3} \right),\quad l=1,\dots,N_j}$;
}
\ENDFOR
\STATE{
$b^w_{k+1}  \gets b^w_k + \Psi u_{k+1} - w_{k+1}$;\\
$b^d_{k+1}  \gets b^d_k + (\nabla^f u_{k+1} -v_{k+1}) -d_{k+1}$;\\
$b^t_{k+1}  \gets b^t_k + \mathcal{E}^bv_{k+1} - t_{k+1}$;\\
$y_{k+1}  \gets y_k + y - Au_{k+1}$;\\
$k \gets k+1$;}
\ENDWHILE
\RETURN Reconstruction $u_{\text{maxIter}}$.
\end{algorithmic}
\end{algorithm*}


We like to comment on two things regarding the algorithm above. First, the weighting matrix $W_j$ depends on the initialization of an $\eps >0$. The choice of $\eps$ is rather empirical and in many cases does not effect the solution. This was already noticed in the beginning of reweighted $\ell^1$ in \cite{CanWakBoy}. This is also the case for our algorithm. As we explained in Section \ref{sec:Section2} the role $\eps$ is essentially just the maximum threshold that our algorithm will do. It provides stability by preventing a division by zero, but the magnitude of the coefficient is mostly determined by the respective analysis coefficient. Let us also mention that for obvious reasons, adapting the $\epsilon$ with respect to the levels cannot achieve the same effect as the proposed choice of the parameters $\lambda_j$.  

Furthermore, we have not incorporated an additional stopping criterion besides a maximum number of iterations. In all our test cases the algorithm shows convergence, see Figure \ref{true}, \ref{GLPU}, \ref{Pepper} and \ref{CT}.

\section{Numerics} \label{sec:Section4}

In this section we will recover numerous signals from their Fourier measurements and  Radon measurements. These are classical problems in applied mathematics where the image of interest is sparse under a multiscale transform such was the wavelet transform or the shearlet transform. One of the most known applications for the recovery of Fourier measurements is magnetic resonance imaging (MRI) where data is collected in the so-called \emph{k-space} which is the Fourier domain, i.e. every point in the k-space can be interpreted as a Fourier coefficient of the object of interest. This is also one of the very first areas where compressed sensing has had a great impact on, see, for instance \cite{LusDonPau}. Our implementation is, however, not restricted to these two types of measurements or sparsifying transforms and other cases can be directly tested.


The numerical performance of our algorithm will be judged by the following three criteria that will be the main points of our analysis:
\begin{itemize}
\item[(N1)] Quality of the reconstruction,
\item[(N2)] Convergence of the algorithm,
\item[(N3)] Stability towards change of signals.
\end{itemize}
 For (N1) we compare our algorithm with different existing and established methods that are known to perform well in the recovery problem from Fourier measurements these methods are shown in Table \ref{Abbrev}. The results are then shown in Section \ref{Sec:Compare}. Moreover, in some cases we will only use a Haarwavelet-based wavelet transform and in some other cases the shearlet transform. This is because we are not interested in promoting a particular transform, but rather show the effectiveness of our algorithm. Therefore we have always chosen the transform that performs best on the specific data set. Only in Figure \ref{Compare} we will show both methods in a direct comparison.
 
 For (N2) we have chosen two quality measurements. First, the relative error which is computed by the formula
 \begin{align*}
    \text{RE} = \frac{\| u^{\text{ref}} - u^{\text{rec}} \|_2}{\| u^{\text{ref}}\|_2},
 \end{align*}
 where $u^{\text{ref}}$ is the reference image and $u^{\text{rec}}$ the reconstructed image, both in a vectorized form. Second, we use the structured similarity index that as introduced in \cite{SSIM}. 
 
 The stability (N3) is verified by the fact that we have chosen the same parameters for each multiscale transform across all experiments for each problem. Although an extensive tuning of all parameters for different images might yield superior results we have chosen not to do so. The reason behind is two-fold: First, our algorithm already performance very well with a fixed choice of parameters for all different images used in this section. Second, iterative reweighting combined with the proposed multilevel weighting strategy already suggests the level of thresholds for all coefficients and should therefore be less sensitive to the choice of additional parameters. 
 
\subsection{Recovery from Fourier measurements}

In this section we give a description of which multiscale transforms we used precisely and the parameters.

For \emph{wavelets} we have used the undecimated 2D wavelet transform of the Spot package available at
\begin{center}
 \texttt{http://www.cs.ubc.ca/labs/scl/spot/}
\end{center}
Otherwise differently stated we have used a 4 scale Daubechies-2 wavelet system.  
The chosen parameters for the proposed algorithm (WIRL1 + TGV) are
\begin{itemize}
\item[-] $\mu$: [6e2 1e1 2e1],
\item[-] $\alpha$: [1 2],
\item[-] $\beta$: 1e4,
\item[-] $\eps$: 1e-4.
\end{itemize}

For \emph{shearlets} we have used the shearlet transform available from
\begin{center}
 \texttt{http://www.shearlab.org/}
\end{center}
The discrete shearlet system is generated by using 4 scales and $[1 \, 1 \,  2 \, 2]$ for the directional parameters. The chosen parameters for the proposed algorithm (SIRL1 + TGV) are
\begin{itemize}
\item[-] $\mu$: [5e3 1e1 2e1],
\item[-] $\alpha$: [1 1],
\item[-] $\beta$: 1e5,
\item[-] $\eps$: 1e-5.
\end{itemize}
For all experiments we chose 2 number of block Gauss-Seidel iterations  and performed 4 inner iterations before updating  $y_k$.  All experiments were conducted in MATLAB R2015b with an Intel i3 CPU with 8GB memory. 

\subsection{Comparison with other methods}\label{Sec:Compare}

\begin{table*}
\centering
\begin{tabularx}{1\textwidth}{|X|X|}
\hline
Fourier inverse & Fourier inversion of data\\
\hline
RecPF & Total variation and wavelet regularization from \cite{RecPF}\\
\hline
FFST+TGV & Shearlets with TGV from \cite{GuoQinYin}\\
\hline
Co-IRL1 & Composite iterative reweighting from \cite{AhmSch}\\
\hline
PANO & Patch-based nonlocal operator, \cite{QuHouLamGuoZhongChen}\\
\hline
TV & Total variation with our algorithm \\
\hline
TGV &  Total generalized variation with our algorithm\\
\hline
WL1 & Wavelets without reweighting\\
\hline
WIRL1 & Wavelets with proposed reweighting\\
\hline
WIRL1+TGV & Wavelets with proposed reweighting and TGV\\
\hline
SL1 & Shearlets without reweighting\\
\hline
SR + TGV & Shearlets with standard reweighting and TGV\\
\hline
SIRL1 & Shearlets with proposed reweighting\\
\hline
SIRL1+TGV & Shearlets with proposed reweighting and TGV\\
\hline
\end{tabularx}
\caption{Table for abbreviations}\label{Abbrev}
\end{table*}

Our first experiment shows the recovery from a $256 \times 256$ rose image available from the open source framework
\begin{center}
 \texttt{http://aforgenet.com/framework/}
\end{center}
We took 30 radial lines through the k-space origin ($\approx 12,2\%$ of Fourier data) representing the measurement data.  We then compared our results obtained by our algorithm using wavelets as well as shearlets.  For both results we have used the proposed iterative multi-level reweighting and a generalized total variation regularizer. The results are compared to RecPF by Yang et al. \cite{RecPF}, Co-IRL1 by Ahmad and Schniter \cite{AhmSch}, PANO by Qu et al. \cite{QuHouLamGuoZhongChen}, and FFST+TGV by Guo et al. \cite{GuoQinYin}. 
In order to make the experiments comparable we have used the same scales and number of directions in \cite{GuoQinYin}. Furthermore, for Co-IRL1 we have used two redundant Daubechies wavelet dictionaries with the same number of scales. Moreover, one dictionary consists of Haar wavelets and the second one of Daubechies 2 wavelets. 

\begin{figure*}
\hspace*{-.9cm}
\begin{tabular}{cccc}
\hspace*{0.4cm}\includegraphics[width=0.235\textwidth]{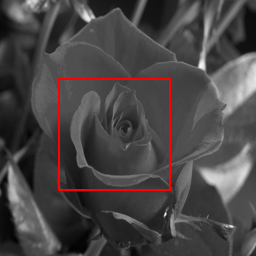} & \hspace*{-0.4cm}
\includegraphics[width=0.235\textwidth]{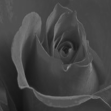}
&
\includegraphics[width=0.235\textwidth]{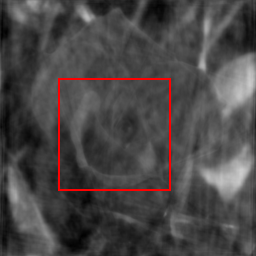}
&
\hspace*{-0.4cm}
\includegraphics[width=0.235\textwidth]{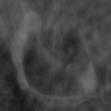}
\\
\multicolumn{2}{c}{\footnotesize{Original}}& \multicolumn{2}{c}
{\footnotesize{Fourier inverse: \quad RE = 0.167, SSIM = 0.739 }}\\[1ex]
 \hspace*{0.4cm}\includegraphics[width=0.235\textwidth]{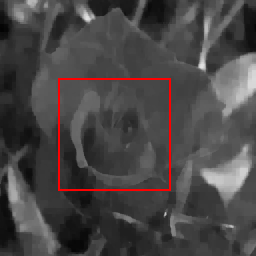} & \hspace*{-0.4cm}
\includegraphics[width=0.235\textwidth]{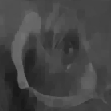}
&
\includegraphics[width=0.235\textwidth]{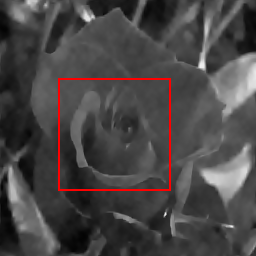}
&
\hspace*{-0.4cm}
\includegraphics[width=0.235\textwidth]{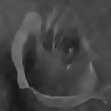}
\\
\multicolumn{2}{c}{\footnotesize{RecPF: \quad RE = 0.105, SSIM = 0.841}} & \multicolumn{2}{c}{\footnotesize{Co-IRL1: \quad RE = 0.076, SSIM = 0.898}}\\[1ex]
\hspace*{0.4cm}\includegraphics[width=0.235\textwidth]{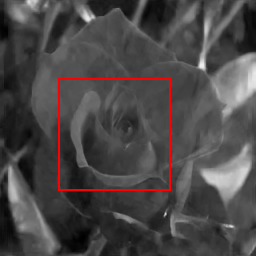}
&
\hspace*{-0.4cm}
\includegraphics[width=0.235\textwidth]{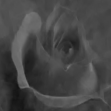}
&
\includegraphics[width=0.235\textwidth]{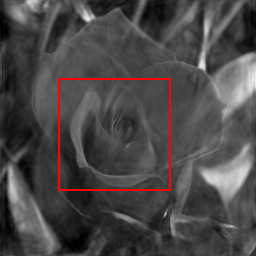}
&
\hspace*{-0.4cm}
\includegraphics[width=0.235\textwidth]{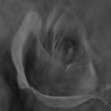}
\\
\multicolumn{2}{c}{\footnotesize{PANO: \quad RE = 0.086, SSIM = 0.881 }} & \multicolumn{2}{c}{\footnotesize{FFST+TGV: \quad RE = 0.095, SSIM = 0.855}} \\[1ex]
\hspace*{0.4cm}\includegraphics[width=0.235\textwidth]{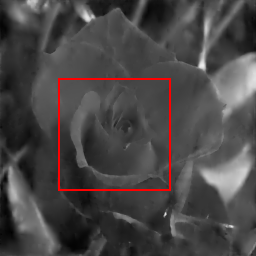}
&
\hspace*{-0.4cm}
\includegraphics[width=0.235\textwidth]{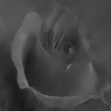}
&
\includegraphics[width=0.235\textwidth]{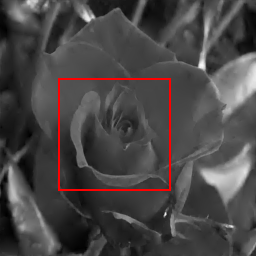}
&
\hspace*{-0.4cm}
\includegraphics[width=0.235\textwidth]{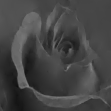}
\\
 \multicolumn{2}{c}{\footnotesize{TGV: \quad RE = 0.113, SSIM = 0.874}} &  \multicolumn{2}{c}{\footnotesize{WIRL1+TGV: \quad RE = 0.058, SSIM = 0.936}}
\end{tabular}
\caption{Different reconstructions from 30 radial lines ($\approx 12,2\%$) through the k-space origin with relative error and structured similarity index. 50 iterations are used for the reconstruction. See Table \ref{Abbrev} for used abbreviations.} \label{Compare}
\end{figure*}

The final result comparing all methods can be found in Figure \ref{Compare}. It can be observed that the recovery obtained
by the proposed method shows the least amount of artifacts and overall the highest quality of reconstruction.

We will next extensively study our algorithm in terms of the effectiveness and stability.

\subsection{Convergence, signal independence, and the effect of reweighting}\label{Sec:Convergence}

\begin{figure*}
\footnotesize
\begin{tabular}{ccccc}
\includegraphics[width=0.17\textwidth]{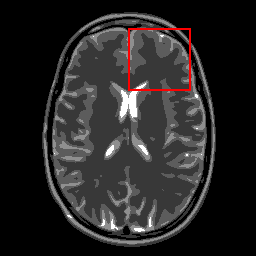}
&\includegraphics[width=0.17\textwidth]{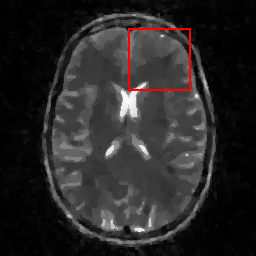}
&\includegraphics[width=0.17\textwidth]{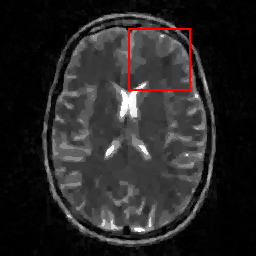}
&\includegraphics[width=0.17\textwidth]{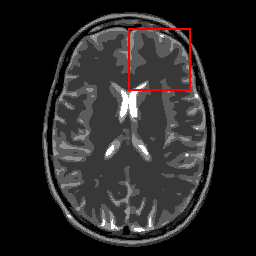}
&\includegraphics[width=0.17\textwidth]{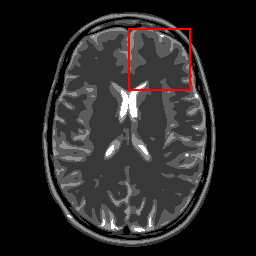}
\\
\footnotesize{Original} & \footnotesize{TV} & \footnotesize{WL1} & \footnotesize{WIRL1} & \footnotesize{WIRL1+TGV} \\
\includegraphics[width=0.17\textwidth]{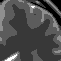}
&\includegraphics[width=0.17\textwidth]{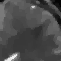}
&\includegraphics[width=0.17\textwidth]{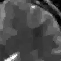}
&\includegraphics[width=0.17\textwidth]{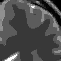}
&\includegraphics[width=0.17\textwidth]{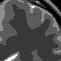}
\\[2ex]
\multicolumn{5}{c}{\includegraphics[width=0.4\textwidth]{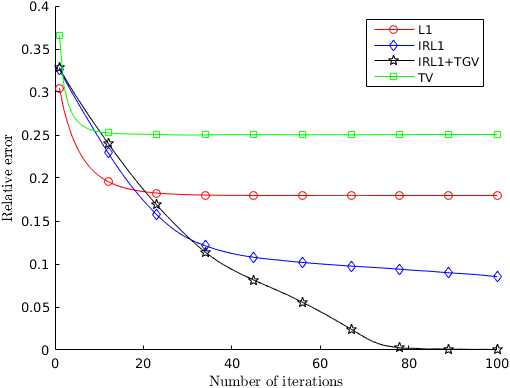}
\includegraphics[width=0.4\textwidth]{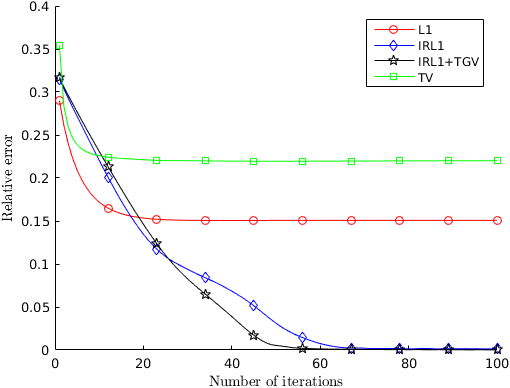}
}
\end{tabular}
\caption{Different reconstructions from 24 radial lines through the k-space origin with relative error and structured similarity index. The lower left graphics corresponds to 21 radial lines ($\approx 8.73 \%$) and the lower right to 24 radial lines ($\approx 9.83 \%$). 100 iterations are used for the reconstruction. See Table \ref{Abbrev} for used abbreviations.} \label{GLPU}
\end{figure*}

In this section we analyze (N2) for our algorithm. We do this by considering two images, one that is well suited for wavelets and the other one where shearlets perform better. We start with a $256 \times 256$ phantom that was designed by Guerquin-Kern et al. in \cite{GLPU} for MRI studies. As this image is piecewise constant we have chosen a 4 scale wavelet transform generated by Haar wavelets. We reconstructed the image using our algorithm for TV only, WL1, WIRL1 and WIRL1+TGV. Moreover, as this image is very compressible in a Haar wavelet basis, the recovery allows a much lower sampling rate. In fact, we have only used 21 radial lines which corresponds to 8.73\%. It is interesting to mention that the exact solution is returned after almost 80 iterations when WIRL1+TGV. Note that, for 24 lines ($\approx 9.83 \%$) wavelets with the proposed iterative reweighting step (WIRL1) will eventually also return the exact solution, see Figure \ref{GLPU}. Also note that the TV reconstruction performs worst, although 
this image should be well suited for TV. Our explanation is that at these extremely low sampling rates a highly redundant transformation is needed to still guarantee recovery.

\begin{figure*}
\begin{tabular}{cccc}
\includegraphics[width=0.22\textwidth]{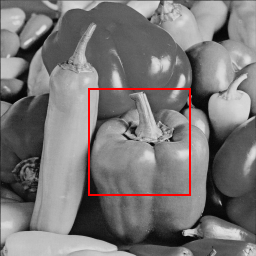}
&\includegraphics[width=0.22\textwidth]{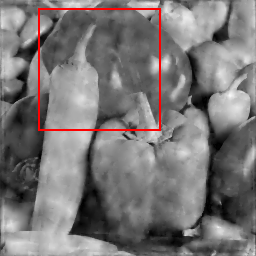}
&\includegraphics[width=0.22\textwidth]{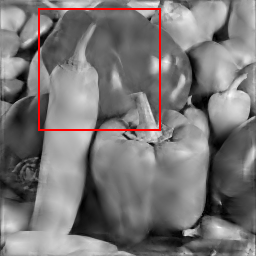}
&\includegraphics[width=0.22\textwidth]{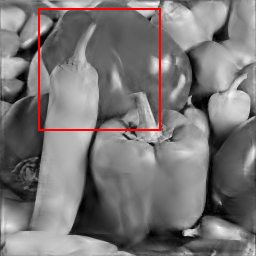}
\\
\footnotesize{Original} & \footnotesize{SL1} & \footnotesize{SR+TGV} & \footnotesize{SIRL1+TGV}
\\
\includegraphics[width=0.22\textwidth]{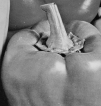}
&\includegraphics[width=0.22\textwidth]{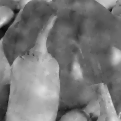}
&\includegraphics[width=0.22\textwidth]{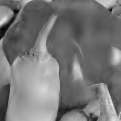}
&\includegraphics[width=0.22\textwidth]{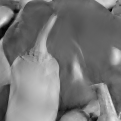}
\\[2ex]
\multicolumn{4}{c}{\includegraphics[width=0.4\textwidth]{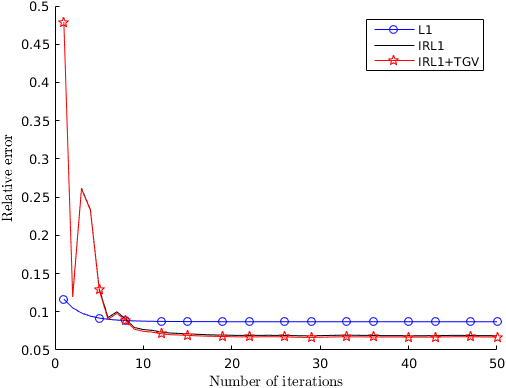}
\includegraphics[width=0.4\textwidth]{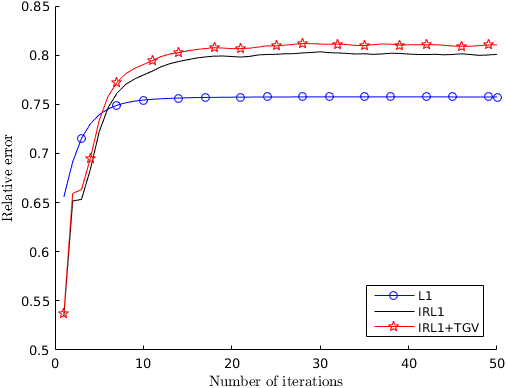}
}
\end{tabular}
\caption{Different reconstruction from 30 radial lines ($\approx 12,2\%$) through the k-space origin with relative error and structured similarity index. See Table \ref{Abbrev} for used abbreviations. 50 iterations are used for the reconstruction.} \label{Pepper}
\end{figure*}

Our third numerical example for Fourier measurements concerns the $256 \times 256$ pepper image, see Figure \ref{Pepper}. It has many more structures than the previously considered GLPU phantom. More importantly, it does not consist of piecewise constant areas. This image is particularly well suited for shearlets and thus we have chosen the shearlet transform  with four scales. The result for a fixed threshold (SL1), i.e. without the proposed reweighting is significantly worse than the one obtained by the multilevel iterative reweighting method (SIRL1). In this case adding total generalized variation as an additional regularizer does not improve the image quality much. Indeed, the improvement of adding TGV as a regularizer depends strongly on the image. It improves the quality if more piecewise constant areas are present in the image (especially if the background is constant as this is typically the case for MRI images), see for instance Figure \ref{GLPU}.

\subsection{Timings}

Although our code is not optimized at all, reconstructions are obtained in reasonable time. Below we display the timings that our algorithm needs in order to compute a reconstruction. The timings are recorded from the experiment shown in Figure \ref{Compare} (image size: $256\times 256$). As our algorithm is also capable of computing reconstructions using only a TV- or TGV-regularizer we additionally recorded these timings.

\begin{table}[H]
\centering
\begin{tabular}{|c|c|}
\hline
 Method & Time in seconds\\
\hline
TV & 16.97\\
\hline
TGV &  47.01\\
\hline
WL1 & 254.77 \\
\hline
WIRL1 & 310.78\\
\hline
WIRL1+TGV & 350.61 \\
\hline
SL1 & 286.57\\
\hline
SIRL1 & 419.20\\
\hline
SIRL1+TGV & 479.74\\
\hline
\end{tabular}
\caption{Table of timings.}\label{Timings}
\end{table}

Note that TGV is more time consuming than TV, mostly due to the larger system \eqref{diagsys}. Moreover, as we already mentioned, the wavelet transform as well as the shearlet transform used in this paper are redundant transforms, i.e. the number of coefficients is in this case significantly larger than the number of pixels available. This explains the much slower performance in comparison to the variational methods.

We also like to mention that the additional cost of reweighting (mostly due to a rather unsophisticated implementation) seems manageable, especially given the benefit of quality that one gets in the recovered images.

\subsection{Recovery from Radon measurements} \label{sec:Radon}
In order to demonstrate the versatility of our proposed method, we furthermore consider an example of a recovery from few Radon measurements. This is another classical linear inverse problem, which forms the underlying model for \emph{computed tomography}. As a test signal we again use the brain phantom of \cite{GLPU}. For the measurement system we use  \emph{fan-beam} projections with 45 equally spaced rotation angles. For the reconstruction of the piecewise constant test signal we make use of the same 4 scale redundant Daubechies-2 wavelet system as above, together with the following parameter setup:
\begin{itemize}
\item[-] $\mu$: [1e3 1e2 2e3],
\item[-] $\alpha$: [1e-3 2e-3],
\item[-] $\beta$: 1e1,
\item[-] $\eps$: 1e-6.
\end{itemize}
Since the matrix $A^*A$ is not diagonalized by $\mathcal{F}$ anymore, we follow the the steps described in Section \ref{sec:Section3} for solving \eqref{diagsys}. The resulting linear system is solved iteratively using 75 \texttt{pcg} iterations. 

\begin{figure*}

\begin{tabular}{cccc}
\hspace*{0.4cm}\includegraphics[width=0.215\textwidth]{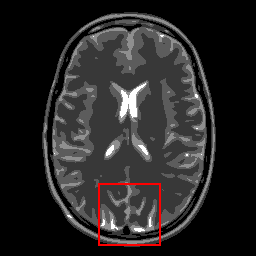} & \hspace*{-0.4cm}
\includegraphics[width=0.215\textwidth]{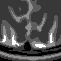}
&
\includegraphics[width=0.215\textwidth]{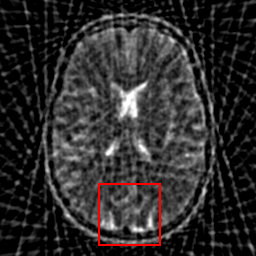}
&
\hspace*{-0.4cm}
\includegraphics[width=0.215\textwidth]{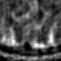}
\\
\multicolumn{2}{c}{\footnotesize{Original}}& \multicolumn{2}{c}
{\footnotesize{Filtered backprojection: \quad RE = 0.507, SSIM = 0.277 }}\\[1ex]
 \hspace*{0.4cm}\includegraphics[width=0.215\textwidth]{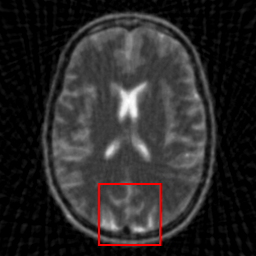} & \hspace*{-0.4cm}
\includegraphics[width=0.215\textwidth]{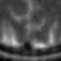}
&
\includegraphics[width=0.215\textwidth]{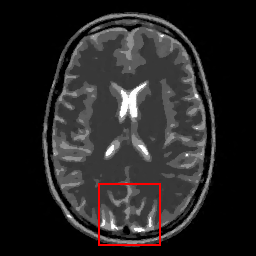}
&
\hspace*{-0.4cm}
\includegraphics[width=0.215\textwidth]{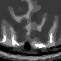}
\\
\multicolumn{2}{c}{\footnotesize{Tikhonov: \quad RE = 0.264, SSIM = 0.479}} & \multicolumn{2}{c}{\footnotesize{TV: \quad RE = 0.094, SSIM = 0.966}}\\[1ex]
\hspace*{0.4cm}\includegraphics[width=0.215\textwidth]{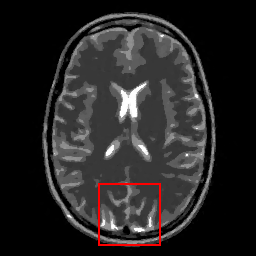}
&
\hspace*{-0.4cm}
\includegraphics[width=0.215\textwidth]{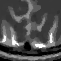}
&
\includegraphics[width=0.215\textwidth]{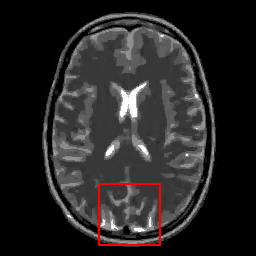}
&
\hspace*{-0.4cm}
\includegraphics[width=0.215\textwidth]{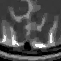}
\\
\multicolumn{2}{c}{\footnotesize{TGV: \quad RE = 0.094, SSIM = 0.967}} & \multicolumn{2}{c}{\footnotesize{WL1: \quad RE = 0.101, SSIM = 0.956}} \\[1ex]
\hspace*{0.4cm}\includegraphics[width=0.215\textwidth]{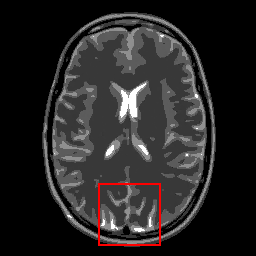}
&
\hspace*{-0.4cm}
\includegraphics[width=0.215\textwidth]{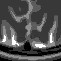}
&
\includegraphics[width=0.215\textwidth]{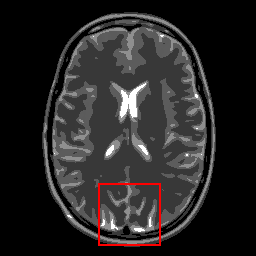}
&
\hspace*{-0.4cm}
\includegraphics[width=0.215\textwidth]{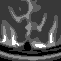}
\\
 \multicolumn{2}{c}{\footnotesize{WIRL1: \quad RE = 0.058 SSIM = 0.985}} &  \multicolumn{2}{c}{\footnotesize{WIRL1+TGV: \quad RE = 0.009, SSIM = 0.999}}
 \\[2ex]
\multicolumn{4}{c}{\includegraphics[width=0.42\textwidth]{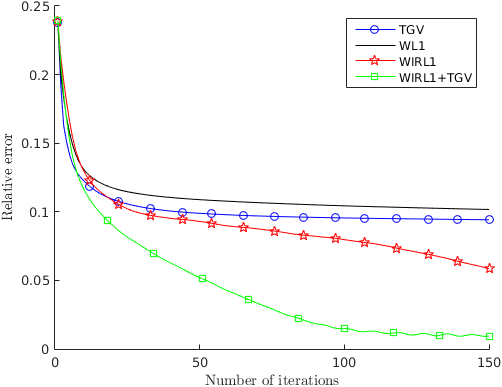}
}
\end{tabular}
\caption{Different reconstructions from 45 angles fan-beam projections with relative error and structured similarity index. 150 iterations are used for the last five methods. See Table \ref{Abbrev} for used abbreviations.} \label{CT}
\end{figure*}

We compare the result with WL1, WIRL1, TV, TGV, Tikhonov regularization and filtered backprojection using a \texttt{`Hann'} filter as provided in Matlab. 
The results are displayed in Figure \ref{CT} and the timings (for 150 iterations) are displayed in Table \ref{TimingsCT}. It can be observed that the convergence is slower as in the Fourier example. The fact that the timings are much higher than before is mostly due to the fact that \eqref{diagsys} is not diagonalizable. Although the Radon-Wavelet measurement matrix is known to be coherent \cite{Sampling,GuaCzaAroLea} (and thus classical CS theorems do not apply), we find it quite remarkable that the measurements contain enough information for an (almost) perfect recovery. 

\begin{table}
\centering
\begin{tabular}{|c|c|}
\hline
 Method & Time in seconds\\
\hline
Filtered backprojection & 0.5\\
\hline
Tikhonov &  0.67\\
\hline
TV & 718 \\
\hline
TGV & 1208\\
\hline
WL1 & 533 \\
\hline
WIRL1 & 703\\
\hline
WIRL1+TGV & 1117\\
\hline
\end{tabular}
\caption{Table of timings for Radon measurements.}\label{TimingsCT}
\end{table}

\section{Conclusion}\label{sec:Section5}


In this paper we presented a novel split Bregman based algorithm that incorporates an iteratively reweighted shrinkage step that is specifically adapted to the structure of multiscale coefficients in order to enhance the quality of image reconstruction. We presented an extensive study of the algorithm thereby focusing on the quality of the image that is to be reconstructed as well as the stability of the parameters towards the change of different signals. Although our focus in the numerical experiments was on the reconstruction problem from Fourier measurements, our algorithm can also be applied to other common problems such as denoising and inpainting. 

Even though all experiments where based on a 2D multiscale transform our findings can be carried over to higher dimensions. This is demonstrated in case of real-world 3D MR data with noise in a practical paper \cite{MaMaeFunSchKutSchKol} using 3D shearlets, \cite{3DShearlets}.

\section*{Acknowledgements}
This work was done when J. Ma was affiliated with the Technische Universit\"{a}t Berlin, he acknowledges support from the DFG Collaborative Research Center TRR 109 “Discretization in Geometry and Dynamics”. JM and MM acknowledge support from the Berlin Mathematical School. The authors would further like to thank Gitta Kutyniok for helpful discussions. Moreover, the authors would like to thank Professor Rizwan Ahmad and Professor Philip Schniter for providing a code for the Co-IRL1 algorithm.

\bibliographystyle{abbrv}
\bibliography{MIR}

\end{document}